\documentclass[12pt,a4paper]{article}

\usepackage{theorem,enumerate}
\usepackage{amsmath,latexsym,amssymb,amsfonts}
\usepackage{graphicx}

\newcounter{memory}

\theorembodyfont{\normalfont\slshape}
\newtheorem{thm}{Theorem}

\newtheorem{Thm}{Theorem}

\newtheorem{cor}[thm]{Corollary}

\newtheorem{lem}{Lemma}

\theoremstyle{remark}

{\theorembodyfont{\upshape}
\newtheorem{rem}[thm]{Remark}}

\newcommand{\qed}{{$\quad\square$\vs{3.6}}}

\newcommand{\vs}[1]{\vspace*{#1 mm}}

\newcommand{\ora}{\overrightarrow}
\newcommand{\ola}{\overleftarrow}

\bfseries\normalfont
\begin{document}

\title{On degree sum conditions\\ for 2-factors with a prescribed number of cycles
}

\author{Shuya Chiba\thanks{E-mail address: \texttt{schiba@kumamoto-u.ac.jp}; This work was supported by JSPS KAKENHI Grant Number 17K05347.}
\vspace{+8pt}
 \\
\small\textsl{Applied Mathematics}\\ 
\small\textsl{Faculty of Advanced Science and Technology}\\ 
\small\textsl{Kumamoto University}\\ 
\vspace{+6pt}
\small\textsl{2-39-1 Kurokami, Kumamoto 860-8555, Japan}
}

\date{}
\maketitle

\vspace{-24pt}
\begin{abstract}
For a vertex subset $X$ of a graph $G$, 
let $\Delta_{t}(X)$ be the maximum value of the degree sums of the subsets of $X$ of size $t$. 
In this paper, 
we prove the following result: 
Let $k$ be a positive integer, and let $G$ be an $m$-connected graph 
of order $n \ge 5k - 2$. 
If $\Delta_{2}(X) \ge n$ for every independent set $X$ of size $\lceil m/k \rceil+1$ in $G$, 
then $G$ has a 2-factor with exactly $k$ cycles. 
This is a common generalization of the results obtained by 
Brandt et al.~[Degree conditions for 2-factors, J.~Graph Theory~24 (1997) 165--173] 
and 
Yamashita~[On degree sum conditions for long cycles and cycles through specified vertices, 
Discrete Math.~308 (2008) 6584--6587], 
respectively.

\end{abstract}

\section{Introduction}
\label{introduction}

In this paper, we consider finite simple graphs, 
which have neither loops nor multiple edges.
For terminology and notation not defined in this paper, we refer the readers to \cite{BMBook2008}. 
The independence number and the connectivity of a graph $G$ are denoted by $\alpha(G)$ and $\kappa(G)$, 
respectively. 
For a vertex $x$ of a graph $G$, 
we denote by $d_{G}(x)$ and $N_{G}(x)$ the degree and the neighborhood of $x$ in $G$. 
Let $\sigma_{m}(G)$ 
be the minimum degree sum of an independent set of $m$ vertices in a graph $G$, 
i.e., 
if $\alpha(G) \ge m$, then 
\begin{align*}
\sigma _{m} (G) 
= 
\min \Big\{ \sum_{x \in X}d_{G}(x) : 
\textup{$X$ is an independent set of $G$ 
with $|X|=m$} \Big\}; 
\end{align*}
otherwise, $\sigma _{m} (G) = +\infty$. 
If the graph $G$ is clear from the context, 
we often omit the graph parameter $G$ in the graph invariant. 
In this paper, ``disjoint'' always means ``vertex-disjoint''.

\medskip
A graph having a \textit{hamilton cycle}, i.e., a cycle containing all the vertices of the graph, 
is said to be \textit{hamiltonian}. 
The hamiltonian problem has long been fundamental in graph theory.
But, it is NP-complete, 
and so 
no easily verifiable necessary and sufficient condition seems to exist.
Therefore, 
many researchers have focused on ``better'' sufficient conditions
for graphs to be hamiltonian (see a survey \cite{LSurvey2013}). 
In particular, the following degree sum condition, due to Ore (1960), is classical and well known.

\begin{Thm}[Ore \cite{Ore1960}]
\label{thm:Ore1960}
Let $G$ be a graph of order $n \ge 3$.
If $\sigma_2 \ge n$, 
then $G$ is hamiltonian.
\end{Thm}

Chv\'{a}tal and Erd\H{o}s (1972) 
discovered 
the relationship between the connectivity, the independence number and the hamiltonicity.

\begin{Thm}[Chv\'{a}tal, Erd\H{o}s \cite{CE1972}]
\label{thm:CE1972}
Let $G$ be a graph of order at least $3$. 
If $\alpha \le \kappa$, 
then $G$ is hamiltonian.
\end{Thm}

Bondy \cite{Bondy1978} pointed out that 
the graph satisfying the Ore condition also satisfies the Chv\'{a}tal-Erd\H{o}s condition, 
that is, 
Theorem~\ref{thm:CE1972} implies Theorem~\ref{thm:Ore1960}. 

By Theorem~\ref{thm:CE1972}, 
we should consider the degree condition 
for the existence of a hamilton cycle in graphs $G$ with $\alpha(G) \ge \kappa (G) + 1$. 
In fact, Bondy (1980) gave the following degree sum condition by extending Theorem~\ref{thm:CE1972}.

\begin{Thm}[Bondy \cite{Bondy1980}]
\label{thm:Bondy1980}
Let $G$ be an $m$-connected graph of order $n \ge 3$. 
If $\sigma_{m + 1} > \frac{1}{2}(m + 1)(n - 1)$, 
then $G$ is hamiltonian.
\end{Thm} 

In 2008, Yamashita \cite{Yamashita2008} introduced a new graph invariant 
and further generalized Theorem~\ref{thm:Bondy1980} as follows. 
For a vertex subset $X$ of a graph $G$ with $|X| \ge t$, we define 
\begin{align*}
\Delta _{t} (X) 
= 
\max \Big\{ \sum_{x \in Y}d_{G}(x) : Y \subseteq X, \ |Y| = t \Big\}. 
\end{align*}
Let $m \ge t$, 
and 
if $\alpha(G) \ge m$, then let 
\begin{align*}
\sigma _{t}^{m} (G) 
= 
\min \Big\{ \Delta _{t} (X)  : 
\textup{$X$ is an independent set of $G$ 
with $|X|=m$} \Big\}; 
\end{align*}
otherwise, $\sigma _{t}^{m} (G) = +\infty$. 
Note that $\sigma_{t}^{m}(G) \ge \frac{t}{m} \cdot \sigma_{m}(G)$.

\begin{Thm}[Yamashita \cite{Yamashita2008}]
\label{thm:Yamashita2008}
Let $G$ be an $m$-connected graph of order $n \ge 3$. 
If $\sigma_2^{m + 1} \ge n$, 
then $G$ is hamiltonian.
\end{Thm} 

This result suggests that 
the degree sum of non-adjacent ``two'' vertices is important for hamilton cycles.

\smallskip
On the other hand, 
it is known that a 2-factor is one of the important generalizations of a hamilton cycle. 
A \textit{2-factor} of a graph is a spanning subgraph in which every component is a cycle, 
and thus a hamilton cycle is a 2-factor with ``exactly 1 cycle''. 
As one of the studies concerning 
the difference between hamilton cycles and 2-factors, 
in this paper, 
we focus on 2-factors with ``exactly $k$ cycles''. 
Similar to the situation for hamilton cycles, 
deciding whether a graph has a 2-factor with $k \ (\ge 2)$ cycles 
is also NP-complete. 
Therefore, 
the sufficient conditions for the existence of such a 2-factor 
also have been extensively studied in graph theory (see a survey \cite{GSurvey2003}). 
In particular, 
the following theorem, due to Brandt, Chen, Faudree, Gould and Lesniak (1997), 
is interesting. 
(In the paper \cite{BCFGL1997}, 
the order condition is not ``$n \ge 4k-1$'' but ``$n \ge 4k$''. 
However, 
by using a theorem of Enomoto \cite{Enomoto1998} and Wang \cite{Wang1999} 
(``every graph $G$ of order at least $3k$ with $\sigma_{2}(G) \ge 4k-1$ 
contains $k$ disjoint cycles'') for the cycles packing problem, 
we can obtain the following. See the proof in \cite[Lemma 1]{BCFGL1997}.)

\begin{Thm}[Brandt et al. \cite{BCFGL1997}]
\label{thm:BCFGL1997}
Let $k$ be a positive integer, 
and let $G$ be a graph of order $n \ge 4k-1$. 
If $\sigma_{2} \ge n$, 
then 
$G$ has a 2-factor with exactly $k$ cycles. 
\end{Thm}

This theorem shows that 
the Ore condition guarantees 
the existence of a hamilton cycle 
but also the existence of a 2-factor with a prescribed number of cycles.

\smallskip
By considering the relation between Theorem~\ref{thm:Ore1960} and Theorem~\ref{thm:BCFGL1997}, 
Chen, Gould, Kawarabayashi, Ota, Saito and Schiermeyer \cite{CGKOSS2007} 
conjectured that the Chv\'{a}tal-Erd\H{o}s condition in Theorem~\ref{thm:CE1972} also guarantees 
the existence of a 2-factor with exactly $k$ cycles (see \cite[Conjecture~1]{CGKOSS2007}). 
Chen et al.~also proved that 
if the order 
of a 2-connected graph $G$ with $\alpha(G) = \alpha \le \kappa(G)$ 
is sufficiently large compared with $k$ and with the Ramsey number $r(\alpha+4, \alpha+1)$, 
then the graph $G$ has a 2-factor with $k$ cycles. 
In \cite{KY2003}, Kaneko and Yoshimoto 
``almost'' solved the above conjecture for $k = 2$ 
(see the comment after Theorem~E in Chen et al.~\cite{CGKOSS2007} for more details). 
Another related result can be found in \cite{CSS2013}. 
But, the above conjecture is still open in general. 
In this sense, there is a big gap between hamilton cycles 
and 2-factors with exactly $k \ (\ge 2)$ cycles.

\medskip
In this paper, 
by combining the techniques of the proof for hamiltonicity 
and the proof for 2-factors with a prescribed number of cycles, 
we give the following 
Yamashita-type condition for 2-factors with $k$ cycles.

\begin{thm}
\label{thm:Yamashita-type for 2-factor with k cycles}
Let $k$ be a positive integer, 
and let $G$ be an $m$-connected graph of order $n \ge 5k - 2$. 
If $\sigma_{2}^{\lceil m/k \rceil + 1} \ge n$, then $G$ has a 2-factor with exactly $k$ cycles. 
\end{thm}

This theorem implies the following:

\begin{rem}
\label{rem:Yamashita-type for 2-factor with k cycles}

\indent
\begin{itemize}
\item 
Theorem~\ref{thm:Yamashita-type for 2-factor with k cycles} is a generalization of Theorem~\ref{thm:Yamashita2008}. 

\item 
Theorem~\ref{thm:Yamashita-type for 2-factor with k cycles} leads to 
the Bondy-type condition: 
If $G$ is an $m$-connected graph of order $n \ge 5k - 2$ with 
$\sigma_{\lceil m / k \rceil + 1}(G) > \frac{1}{2}(\lceil m / k \rceil + 1)(n-1)$, 
then $G$ has a 2-factor with exactly $k$ cycles. 
Therefore, 
Theorem~\ref{thm:Yamashita-type for 2-factor with k cycles} is also a generalization of Theorem~\ref{thm:BCFGL1997} 
for sufficiently large graphs.
(Recall that 
$\sigma_{t}^{m}(G) \ge \frac{t}{m} \cdot \sigma_{m}(G)$ 
and 
$\sigma_{m}(G) \ge \frac{m}{2} \cdot \sigma_{2}(G)$ for $m \ge t \ge 2$.)

\item 
Theorem~\ref{thm:Yamashita-type for 2-factor with k cycles} leads to 
the Chv\'{a}tal-Erd\H{o}s-type condition: 
If $G$ is a graph of order at least $5k - 2$ with $\alpha(G) \le \lceil \kappa(G) / k \rceil$, 
then $G$ has a 2-factor with exactly $k$ cycles. 
\end{itemize}
\end{rem}

The complete bipartite graph $K_{(n-1)/2, (n + 1)/2}$ ($n$ is odd) 
does not contain a 2-factor, 
and hence the degree condition in 
Theorem~\ref{thm:Yamashita-type for 2-factor with k cycles} is best possible in this sense. 
The order condition in Theorem~\ref{thm:Yamashita-type for 2-factor with k cycles} comes from our proof techniques. 
Similar to the situation for the proof of Theorem~\ref{thm:BCFGL1997}, 
we will use the order condition only for the cycles packing problem 
(see Lemma~\ref{lem:cycle packing} and the proof of Theorem~\ref{thm:Yamashita-type for 2-factor with k cycles} in Section~\ref{proof of main}). 
The complete bipartite graph $K_{2k-1, 2k-1}$ shows that $n \ge 4k-1$ is necessary. 
In the last section (Section~\ref{order condition}), 
we note that ``$n \ge 5k-2$'' can be replaced with ``$n \ge 4k-1$'' for the Bondy-type condition (and the Chv\'{a}tal-Erd\H{o}s-type condition) 
in Remark~\ref{rem:Yamashita-type for 2-factor with k cycles}.

\smallskip
Table~\ref{compare} summarizes the conditions mentioned in the above.

\begin{table}[h]
\begin{center}
\begin{tabular}{| c  |  c |  c |}
\hline
&
\footnotesize{hamilton cycle}
&
\footnotesize{2-factor with $k$ cycles}
\\
\hline
\hline
\footnotesize{Ore-type}
&
\footnotesize{$\sigma_{2} \ge n$}
&
\footnotesize{$\sigma_{2} \ge n$}
\\

&
\scriptsize{Theorem~\ref{thm:Ore1960} (Ore)}
&
\scriptsize{Theorem~\ref{thm:BCFGL1997} (Brandt et al.)}
\\
\hline
\footnotesize{Chv\'{a}tal-Erd\H{o}s-type}
&
\footnotesize{$\alpha \le \kappa$}
&
\footnotesize{$\alpha \le \lceil \ \kappa / k \ \rceil$}
\\

&
\scriptsize{Theorem~\ref{thm:CE1972} (Chv\'{a}tal and Erd\H{o}s)}
&
\scriptsize{Remark~\ref{rem:Yamashita-type for 2-factor with k cycles}}
\\
\hline
\footnotesize{Bondy-type}
&
\footnotesize{\hspace{+13pt} $\sigma_{\kappa+1} > \frac{1}{2}(\kappa+1)(n-1)$ \hspace{+13pt}}
&
\footnotesize{$\sigma_{\lceil \kappa / k \rceil+1} > \frac{1}{2}(\lceil \kappa / k \rceil+1)(n-1)$}
\\

&
\scriptsize{Theorem~\ref{thm:Bondy1980} (Bondy)}
&
\scriptsize{Remark~\ref{rem:Yamashita-type for 2-factor with k cycles}}
\\
\hline
\footnotesize{Yamashita-type}
&
\footnotesize{$\sigma_{2}^{\kappa +1} \ge n$}
&
\footnotesize{$\sigma_{2}^{\lceil \kappa / k \rceil+1} \ge n$} 
\\

&
\scriptsize{Theorem~\ref{thm:Yamashita2008} (Yamashita)}
&
\scriptsize{Theorem~\ref{thm:Yamashita-type for 2-factor with k cycles} (Main theorem)}
\\
\hline
\end{tabular}
\caption{Comparison of the degree conditions} 
\label{compare}
\end{center}
\end{table}

\smallskip
To prove Theorem~\ref{thm:Yamashita-type for 2-factor with k cycles}, 
in the next section, 
we extend the concept of insertible vertices which was introduced by Ainouche \cite{Ainouche1992}, 
and we prove Theorem~\ref{thm:Yamashita-type for 2-factor with k cycles} in Section~\ref{proof of main} by using it.

\section{The concept of insertible vertices}
\label{concept of insertible}

In this section, 
we prepare terminology and notations 
and give some lemmas.

Let $G$ be a graph. 
For $v \in V(G)$ and $X \subseteq V(G)$, 
we let $N_{G}(v; X) = N_{G}(v) \cap X$ 
and $d_{G}(v; X) = |N_{G}(v; X)|$. 
For $V, X \subseteq V(G)$, 
let $N_{G}(V; X) = \bigcup_{v \in V}N_{G}(v; X)$. 
For $X \subseteq V(G)$, 
we denote by $G[X]$ the subgraph of $G$ induced by $X$. 
An $(x, y)$-\textit{path} in $G$ is a path from a vertex $x$ to a vertex $y$ in $G$ . 
We write a cycle (or a path) $C$ with a given orientation by $\ora{C}$. 
If there exists no fear of confusion, 
we abbreviate $\ora{C}$ by $C$. 
Let $C$ be an oriented cycle (or path). 
We denote by $\ola{C}$ the cycle $C$ with a reverse orientation. 
For $x \in V(C)$, 
we denote the successor and the predecessor of $x$ on $\ora{C}$ 
by $x^{+}$ and $x^{-}$. 
For $x, y \in V(C)$,
we denote by $x\ora{C}y$
the $(x, y)$-path on $\ora{C}$.
The reverse sequence of $x\ora{C}y$ is denoted by $y\ola{C}x$. 
In the rest of this paper, 
we consider that 
every cycle (path) has a fixed orientation, unless stated otherwise, 
and 
we often identify a subgraph $H$ of $G$ with its vertex set $V(H)$.

\medskip
The following lemma is obtained by using the standard crossing argument, 
and so we omit the proof.

\begin{lem}
\label{lem:standard crossing}
Let $G$ be a graph of order $n$, 
and let $P$ be an $(x, y)$-path of order at least $3$ in $G$.  
If $d_{G}(x) + d_{G}(y) \ge n$, 
then $G$ contains a cycle of order at least $|P|$. 
\end{lem}

In \cite{Ainouche1992}, 
Ainouche introduced the concept of insertible vertices, 
which has been used for the proofs of the results on hamilton cycles. 
In this paper, 
we modify it for 2-factors with $k$ cycles, 
and it also plays a crucial role in our proof. 
Let $G$ be a graph, 
and let $\mathcal{D} = \{D_{1}, \dots, D_{r+s}\}$ $(r+s \ge 1)$ be the set of $r$ cycles and $s$ paths in $G$ 
which are pairwise disjoint. 
For a vertex $x$ in $G - \bigcup_{1 \le p \le r+s}D_{p}$, 
the vertex $x$ is \textit{insertible for $\mathcal{D}$} 
if there is an edge $uv$ in $E(D_{p})$ 
such that $xu, xv \in E(G)$ for some $p$ with $1 \le p \le r+s$. 
In the following lemma, 
``partition'' of a graph means a partition of the vertex set.

\begin{lem}
\label{lem:insertible}
Let $G$ be a graph, 
and let $\mathcal{D} = \{D_{1}, \dots, D_{r+s}\}$ $(r+s \ge 1)$ be the set of $r$ cycles and $s$ paths in $G$ 
which are pairwise disjoint, 
and 
$P$ be a path in $G - \bigcup_{1 \le p \le r+s}D_{p}$. 
If every vertex of $P$ is insertible for $\mathcal{D}$, 
then 
$G \big[ \bigcup_{1 \le p \le r+s}V(D_{p}) \cup V(P) \big]$ 
can be partitioned into $r$ cycles and $s$ paths. 
\end{lem}
\noindent
\textbf{Proof of Lemma~\ref{lem:insertible}.}~By choosing 
the following two vertices $u, v \in V(P)$ and the edge $ww^{+} \in \bigcup_{1 \le p \le r+s}E(D_{p})$ inductively,
we can get the desired partition of $G \big[ \bigcup_{1 \le p \le r+s}V(D_{p}) \cup V(P) \big]$. 
Let 
$u$ be the first vertex along $\ora{P}$, 
and take an edge $ww^{+}$ in $E(D_{i}) \ \big(\subseteq \bigcup_{1 \le p \le r+s}E(D_{p}) \big)$  
such that $uw, uw^{+} \in E(G)$ for some $i$ with $1 \le i \le r+s$ 
(since $u$ is insertible for $\mathcal{D}$, 
we can take such an edge). 
We let $v$ be the last vertex along $\ora{P}$ 
such that 
$vw, vw^{+} \in E(G)$ (may be $u = v$). 
Then, 
we can insert all vertices of $u \ora{P} v$ into $D_{i}$. 
In fact, 
by replacing the edge $ww^{+}$ by the path $wu \ora{P} v w^{+}$, 
we can obtain a spanning subgraph $D_{i}'$ of $G[V(D_{i} \cup u \ora{P} v)]$ 
such that $D_{i}'$ is a cycle if $D_{i}$ is a cycle; 
otherwise, $D_{i}'$ is a path. 
By the choice of $u$ and $v$, 
we have $zw \notin E(G)$ or $zw^{+} \notin E(G)$ 
for each vertex $z$ of $P' := P - u\ora{P}v$, 
and hence 
every vertex of $P'$ 
is insertible for $\mathcal{D'} = \{D_{1}, \dots, D_{i-1}, D_{i}', D_{i + 1}, \dots, D_{r+s}\}$. 
Thus, we can repeat this argument 
for the path $P'$ and the set $\mathcal{D}'$, 
and we get then the desired partition. 
\qed

In the rest of this section, 
we fix the following. 
Let $C_{1}, \dots, C_{k}$ be $k$ disjoint cycles in a graph $G$, 
and let $C^{*} = \bigcup_{1 \le p \le k}C_{p}$. 
Choose $C_{1}, \dots, C_{k}$ so that 
\begin{align*}
|C^{*}| \ \big( = \sum_{1 \le p \le k}|C_{i}| \big) \textup{ is as large as possible.}
\end{align*}
Suppose that $C^{*}$ does not form a 2-factor of $G$. 
Let $H = G - C^{*}$, 
and let $H_{0}$ be a component of $H$ and $x_{0} \in V(H_{0})$. 
Let 
\begin{align*}
\textup{$u_{1}, u_{2}, \dots, u_{l}$ be $l$ distinct vertices in $N_{G}(H_{0}; C_{1})$, where $l \ge 2$.}
\end{align*}
We assume that 
$u_{1}, u_{2}, \dots, u_{l}$ appear in this order on $\ora{C_{1}}$, 
and let $u_{l+1} = u_{1}$. 
Note that by the maximality of $|C^{*}|$, $u_{i}^{+} \neq u_{i+1}$ for $1 \le i \le l$. 
We denote by $\ora{Q_{i}}$ and $\ora{Q_{i, j}}$
a $(u_{i}, x_{0})$-path in $G[V(H_{0}) \cup \{u_{i}\}]$ 
and 
a $(u_{i}, u_{j})$-path passing through a vertex of $H_{0}$ in $G[V(H_{0}) \cup \{u_{i}, u_{j}\}]$, 
respectively.

\begin{lem}
\label{lem:noninsertible vertex}
For $1 \le i \le l$, 
$u_{i}^{+} \ora{C_{1}} u_{i+1}^{-}$ contains a non-insertible vertex for $\{C_{2}, \dots, C_{k}\}$. 
\end{lem}
\noindent
\textbf{Proof of Lemma~\ref{lem:noninsertible vertex}.}~Suppose that 
every vertex of $u_{i}^{+} \ora{C_{1}} u_{i+1}^{-}$ is insertible for $\{C_{2}, \dots,$ $C_{k}\}$. 
Then, 
by Lemma~\ref{lem:insertible}, 
$G \big[ \bigcup_{2 \le p \le k}V(C_{p}) \cup V(u_{i}^{+} \ora{C_{1}} u_{i+1}^{-}) \big]$ 
has a 2-factor with exactly $k-1$ cycles. 
With the cycle $u_{i+1} \ora{C_{1}} u_{i} \ora{Q_{i, i+1}} u_{i+1}$, 
we can get $k$ disjoint cycles in $G$ such that 
the sum of the orders is larger than $|C^{*}|$, a contradiction. 
\qed

For $1 \le i \le l$, 
let $x_{i}$ be the first non-insertible vertex for $\{C_{2}, \dots, C_{k}\}$ 
in $V(u_{i}^{+} \ora{C_{1}} u_{i+1}^{-})$ on $\ora{C_{1}}$, 
i.e., 
every vertex of $u_{i}^{+} \ora{C_{1}} x_{i}^{-}$ 
is insertible for $\{C_{2}, \dots, C_{k}\}$, 
but $x_{i}$ is not insertible 
(Lemma~\ref{lem:noninsertible vertex} guarantees the existence of such a vertex $x_{i}$).

\begin{lem}
\label{lem:crossing}
Let $i, j$ be integers with $1 \le i, j \le l$ and $i \neq j$. 
If $x \in V(u_{i}^{+} \ora{C_{1}} x_{i})$ and $x' \in \{x_{0}, u_{j}^{+}\}$, 
then {\rm (i)} $xx' \notin E(G)$, and {\rm (ii)} $d_{G}(x; H \cup C_{1}) + d_{G}(x'; H \cup C_{1}) \le |H \cup C_{1}| - 1$.  
\end{lem}
\noindent
\textbf{Proof of Lemma~\ref{lem:crossing}.}~Consider the path 
\begin{align*}
\ora{P} 
= 
\left \{
\begin{array}{ll} 
x \ora{C_{1}} u_{i} \ora{Q_{i}} x_{0} & \textup{(if $x' = x_{0}$)} \\[3mm]
x\ora{C_{1}}u_{j} \ola{Q_{i, j}} u_{i} \ola{C_{1}} u_{j}^{+} & \textup{(if $x' = u_{j}^{+}$)}
\end{array}
\right..
\end{align*}
\begin{figure}[h]
\begin{center}
\includegraphics[scale=0.35,clip]{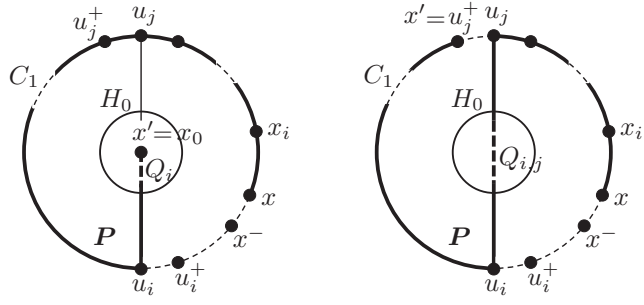}
\caption{The path $P$}
\label{fig:crossing}
\end{center}
\end{figure}
See Figure~\ref{fig:crossing}.  
Then, $P$ is a path in $G[V(H \cup x \ora{C_{1}} u_{i})]$ 
passing through all vertices of $x \ora{C_{1}} u_{i}$ and a vertex of $H_{0}$. 
Recall that every vertex of $u_{i}^{+} \ora{C_{1}} x^{-}$ 
is insertible for $\{C_{2}, \dots, C_{k}\}$, 
and hence $G \big[ \bigcup_{2 \le p \le k}V(C_{p}) \cup V(u_{i}^{+} \ora{C_{1}} x^{-}) \big]$ has a 2-factor with exactly $k-1$ cycles 
(by Lemma~\ref{lem:insertible}). 
Hence, the maximality of $|C^{*}|$ and Lemma~\ref{lem:standard crossing} yield that 
\begin{align*}
xx' \notin E(G)
\textup{ and }
d_{G}(x; H \cup x \ora{C_{1}} u_{i}) + d_{G}(x'; H \cup x \ora{C_{1}} u_{i}) \le |H \cup x \ora{C_{1}} u_{i}| - 1.  
\end{align*}
In particular, 
(i) holds. 
Then, by applying (i) for each vertex in $u_{i}^{+} \ora{C_{1}} x^{-}$ and the vertex $x'$, 
we have $N_{G}(x'; u_{i}^{+} \ora{C_{1}} x^{-}) = \emptyset$. 
Combining this with the above inequality, 
we get, 
\begin{align*}
&\hspace{+16pt}
d_{G}(x; H \cup C_{1}) + d_{G}(x'; H \cup C_{1}) \\
&= d_{G}(x; H \cup x \ora{C_{1}} u_{i}) + d_{G}(x'; H \cup x \ora{C_{1}} u_{i}) + d_{G}(x; u_{i}^{+} \ora{C_{1}} x^{-})\\ 
&\le \big( |H \cup x \ora{C_{1}} u_{i}| - 1 \big) + |u_{i}^{+} \ora{C_{1}} x^{-}| = |H \cup C_{1}| - 1. 
\end{align*}
Thus (ii) also holds. 
\qed

\section{Proof of Theorem~\ref{thm:Yamashita-type for 2-factor with k cycles}}
\label{proof of main}

Before proving Theorem~\ref{thm:Yamashita-type for 2-factor with k cycles}, 
we will give the following lemma for the cycles packing problem.

\begin{lem}
\label{lem:cycle packing}
Let $k, m, n$ and $G$ be the same ones as in Theorem~\ref{thm:Yamashita-type for 2-factor with k cycles}. 
Under the same degree sum condition as Theorem~\ref{thm:Yamashita-type for 2-factor with k cycles}, 
$G$ contains $k$ disjoint cycles. 
\end{lem}
\medskip
\noindent
\textbf{Proof of Lemma~\ref{lem:cycle packing}.}~If $k = 1$, 
then it is easy to check that $G$ contains a cycle. 
If $\lceil m / k \rceil = 1$ or $\lceil m / k \rceil \ge 3$, 
then by a theorem of Enomoto \cite{Enomoto1998}, 
$G$ contains $k$ disjoint cycles 
(note that 
if $\lceil m / k \rceil \ge 3$, then 
$G$ is $(2k+1)$-connected, 
that is, the minimum degree $\delta(G)$ is at least $2k + 1$). 
Thus, we may assume that $k \ge 2$ and $\lceil m / k \rceil = 2$. 
Then, 
we have $\delta(G) \ge m \ge k + 1$ 
and 
$\sigma_{2}^{3}(G) = \sigma_{2}^{\lceil m / k \rceil + 1}(G) \ge n \ge 5k - 2$. 
Note that, 
by the definition of $\sigma_{2}^{3}(G)$ and $\sigma_{3}(G)$, 
$\sigma_{3}(G) \ge \sigma_{2}^{3}(G) + \delta(G)$. 
Note also that $n \ge 5k - 2 \ge 3k + 2 \ge 8$ because $k \ge 2$. 
Hence, 
by a theorem of Fujita et al.~\cite{FMTY2006} 
(``every graph $G$ of order at least $3k+2 \ge 8$ with $\sigma_{3}(G) \ge 6k - 2$ 
contains $k$ disjoint cycles''), 
we can get the desired conclusion. 
\qed

Now we are ready to prove Theorem~\ref{thm:Yamashita-type for 2-factor with k cycles}.

\medskip
\noindent
\textbf{Proof of Theorem~\ref{thm:Yamashita-type for 2-factor with k cycles}.}~Let $G$ 
be an $m$-connected graph of order $n \ge 5k-2$ 
such that $\sigma_{2}^{\lceil m / k \rceil + 1}(G) \ge n$. 
We show that $G$ has a 2-factor with exactly $k$ cycles. 
By Theorem~\ref{thm:BCFGL1997}, we may assume that $\lceil m / k \rceil \ge 2$. 
By Lemma~\ref{lem:cycle packing}, 
$G$ contains $k$ disjoint cycles. 
Let $C_{i}$ ($1 \le i \le k$), $C^{*}$, $H$, $H_{0}$, $x_{0}$ and $u_{i}$ ($1 \le i \le l$) 
be the same graphs and vertices as the ones described in the paragraph preceding Lemma~\ref{lem:noninsertible vertex} in Section~\ref{concept of insertible}. 
In particular, 
we may assume that 
$l = \lceil m/k \rceil$. 
Because, 
since $G$ is $m$-connected, 
it follows that $|N_{G}(H_{0}; C^{*})| \ge m$ 
(note that by the maximality of $|C^{*}|$, $|C^{*}| > m$), 
and hence, 
without loss of generality, 
we may assume that $|N_{G}(H_{0}; C_{1})| \ge \lceil m / k \rceil \ (\ge 2)$.

\medskip
We first consider the set 
\begin{align*}
X = \{x_{0} \} \cup \{ u_{i}^{+} : 1 \le i \le l \}. 
\end{align*}
Then, Lemma~\ref{lem:crossing} implies the following: 
\begin{enumerate}[{\upshape(1)}]
\setcounter{enumi}{\value{memory}}
\item
\label{X is independent}
$X$ is an independent set of size $l + 1$. 

\item
\label{degree sum of x1 and x2 in H and C1}
$d_{G}(x; H \cup C_{1}) + d_{G}(x'; H \cup C_{1}) \le |H \cup C_{1}| - 1$ for $x, x' \in X$ ($x \neq x'$). 

\setcounter{memory}{\value{enumi}}
\end{enumerate}

On the other hand, 
by the maximality of $|C^{*}|$ and Lemma~\ref{lem:insertible}, 
$x_{0}$ is non-insertible for $\{C_{2}, \dots, C_{k}\}$. 
This implies the following: 
\begin{enumerate}[{\upshape(1)}]
\setcounter{enumi}{\value{memory}}

\item
\label{degree of x0 in Cp} 
$d_{G}(x_{0}; C_{p}) \le |C_{p}|/2$ for $2 \le p \le k$, and hence $d_{G}(x_{0}; C^{*}-C_{1}) \le |C^{*}-C_{1}|/2$. 

\setcounter{memory}{\value{enumi}}
\end{enumerate}

\smallskip
Since $\sigma_{2}^{l+1}(G) \ge n$, 
it follows from (\ref{X is independent}) that there exist two distinct vertices $x$ and $x'$ in $X$ 
such that $d_{G}(x) + d_{G}(x') \ge n$. 
Then, by (\ref{degree sum of x1 and x2 in H and C1}), we get 
\begin{align*}
d_{G}(x; C^{*} - C_{1}) + d_{G}(x'; C^{*} - C_{1}) 
\ge 
n - \big( |H \cup C_{1}| - 1 \big) = |C^{*} - C_{1}| + 1. 
\end{align*}
Combining this with 
(\ref{degree of x0 in Cp}) and the definition of $X$, 
we may assume that 
\begin{enumerate}[{\upshape(1)}]
\setcounter{enumi}{\value{memory}}

\item 
\label{degree of u1+ in C*-C1}
$d_{G}(u_{1}^{+}; C^{*} - C_{1}) > |C^{*} - C_{1}|/2$. 

\setcounter{memory}{\value{enumi}}
\end{enumerate}

\medskip
Next, let $x_{1}$ be the first non-insertible vertex for $\{C_{2}, \dots, C_{k}\}$ 
in the path $u_{1}^{-} \ola{C_{1}} u_{l}^{+}$ on $\ola{C_{1}}$ 
(we can take such a vertex by Lemma~\ref{lem:noninsertible vertex} and the symmetry of $\ora{C_{1}}$ and $\ola{C_{1}}$), 
and 
we consider the set 
\begin{align*}
Y = \{x_{0}, x_{1} \} \cup \{ u_{i}^{-} : 2 \le i \le l \}. 
\end{align*}
Then, 
by the symmetry of $\ora{C_{1}}$ and $\ola{C_{1}}$, Lemma~\ref{lem:crossing}, 
and since $x_{1}$ is non-insertible for $\{C_{2}, \dots, C_{k}\}$, 
we have the following: 
\begin{enumerate}[{\upshape(1)}]
\setcounter{enumi}{\value{memory}}

\item
\label{Y is independent} 
$Y$ is an independent set of size $l + 1$. 

\item
\label{degree sum of y1 and y2 in H and C1}
$d_{G}(y; H \cup C_{1}) + d_{G}(y'; H \cup C_{1}) \le |H \cup C_{1}| - 1$ for $y, y' \in Y$ ($y \neq y'$). 

\item 
\label{degree of x1 in Cp} 
$d_{G}(x_{1}; C_{p}) \le |C_{p}|/2$ for $2 \le p \le k$, 
and hence $d_{G}(x_{1}; C^{*}-C_{1}) \le |C^{*}-C_{1}|/2$. 

\setcounter{memory}{\value{enumi}}
\end{enumerate}

\smallskip
Since $\sigma_{2}^{l+1}(G) \ge n$, 
it follows from (\ref{Y is independent}) that there exist two distinct vertices $y$ and $y'$ in $Y$ 
such that $d_{G}(y) + d_{G}(y') \ge n$. 
Then, by (\ref{degree sum of y1 and y2 in H and C1}), we get 
\begin{align*}
d_{G}(y; C^{*} - C_{1}) + d_{G}(y'; C^{*} - C_{1}) 
\ge 
n - \big( |H \cup C_{1}| - 1 \big) = |C^{*} - C_{1}| + 1. 
\end{align*}
Combining this with 
(\ref{degree of x0 in Cp}), (\ref{degree of x1 in Cp}) and the definition of $Y$, 
we have the following: 
\begin{enumerate}[{\upshape(1)}]
\setcounter{enumi}{\value{memory}}

\item 
\label{degree of ui- in C*-C1}
$d_{G}(u_{i}^{-}; C^{*} - C_{1}) > |C^{*} - C_{1}|/2$ for some $i$ with $2 \le i \le l$. 

\setcounter{memory}{\value{enumi}}
\end{enumerate}

\begin{figure}[h]
\begin{center}
\includegraphics[scale=0.35,clip]{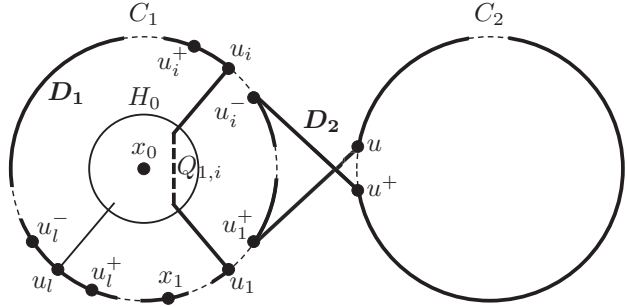}
\caption{The cycles $D_{1}$ and $D_{2}$}
\label{fig:insertiblepath}
\end{center}
\end{figure}

\smallskip
By (\ref{degree of u1+ in C*-C1}) and (\ref{degree of ui- in C*-C1}), 
we have 
\begin{align*}
d_{G}(u_{1}^{+}; C^{*} - C_{1}) + d_{G}(u_{i}^{-}; C^{*} - C_{1}) > |C^{*} - C_{1}| 
= \sum_{2 \le p \le k}|C_{p}|.
\end{align*} 
Hence, there exists a cycle $C_{p}$ ($2 \le p \le k$), 
say $p = 2$, 
such that 
\begin{align*}
d_{G}(u_{1}^{+}; C_{2}) + d_{G}(u_{i}^{-}; C_{2}) \ge |C_{2}| + 1. 
\end{align*}
This implies that 
there exists an edge $uv$ in $E(C_{2})$ such that 
$u_{1}^{+}u, u_{i}^{-}v \in E(G)$. 
By changing the orientation of $C_{2}$ if necessary, 
we may assume that 
$u^{+} = v$. 
Note that $i \ge 2$, 
and consider two cycles 
\begin{align*}
D_{1} = u_{i} \ora{C_{1}} u_{1} \ora{Q_{1, i}} u_{i} 
\textup{ and }
D_{2} = u_{1}^{+} \ora{C_{1}} u_{i}^{-} u^{+} \ora{C_{2}} u u_{1}^{+} \textup{ (see Figure~\ref{fig:insertiblepath})}. 
\end{align*}
Then, 
$D_{1}, D_{2}, C_{3}, \dots, C_{k}$ 
are $k$ disjoint cycles 
such that the sum of the orders is larger than $|C^{*}|$, a contradiction.

This completes the proof of Theorem~\ref{thm:Yamashita-type for 2-factor with k cycles}. 
\qed

\section{Notes on the order condition}
\label{order condition}

As shown in the argument of the previous section, 
in the proof of Theorem~\ref{thm:Yamashita-type for 2-factor with k cycles}, 
the order condition ``$n \ge 5k-2$'' is required only to show the existence of $k$ disjoint cycles in a graph $G$ 
(recall that the order condition in Theorem~\ref{thm:BCFGL1997} is also). 
Therefore, 
the proof of Theorem~\ref{thm:Yamashita-type for 2-factor with k cycles} actually implies the following.

\begin{thm}
\label{thm:Yamashita-type for 2-factor with k cycles (Partition)}
Let $k$ be a positive integer, 
and let $G$ be an $m$-connected graph of order $n$. 
Suppose that $G$ contains $k$ disjoint cycles. 
If $\sigma_{2}^{\lceil m/k \rceil + 1} \ge n$, then $G$ has a 2-factor with exactly $k$ cycles. 
\end{thm}

From this theorem, 
if we can obtain better results on the cycles packing problem, 
then the order conditions in Theorem~\ref{thm:Yamashita-type for 2-factor with k cycles} 
and Remark~\ref{rem:Yamashita-type for 2-factor with k cycles} can be improved. 
In fact, 
by using the result of Kierstead, Kostochka and Yeager (2017) 
and modifying the proof of Lemma~\ref{lem:cycle packing}, 
we can obtain a sharp order condition for the result in Remark~\ref{rem:Yamashita-type for 2-factor with k cycles} 
(see Corollary~\ref{cor:Bondy-type for 2-factor with k cycles}).

\begin{Thm}[Kierstead et al. \cite{KKY2017}]
\label{Thm:KKY2017}
Let $k$ be an integer with $k \ge 2$, 
and let $G$ be a graph of order $n \ge 3k$ with $\delta(G) \ge 2k-1$. 
Then $G$ contains $k$ disjoint cycles if and only if 
{\rm (i)} $\alpha(G) \le n - 2k$, 
and 
{\rm (ii)} if $k$ is odd and $n =3k$, then $G \not \cong 2K_{k} \vee \overline{K_{k}}$ 
and if $k=2$, then $G$ is not a wheel.
\end{Thm}

\begin{lem}
\label{lem:cycle packing II}
Let $k$ be a positive integer, 
and let $G$ be an $m$-connected graph of order $n \ge 4k-1$. 
If $\sigma_{\lceil m/k \rceil + 1}(G) > \frac{1}{2}(\lceil m/k \rceil + 1) (n-1)$, 
then $G$ contains $k$ disjoint cycles. 
\end{lem}
\medskip
\noindent
\textbf{Proof of Lemma~\ref{lem:cycle packing II}.}~By a similar argument 
as in the proof of Lemma~\ref{lem:cycle packing}, 
we have the following: 
If $k = 1$, 
then we can easily find a cycle; 
If $\lceil m / k \rceil = 1$ or $\lceil m / k \rceil \ge 3$, 
then by a theorem of Enomoto \cite{Enomoto1998}, 
$G$ contains $k$ disjoint cycles; 
If $\lceil m / k \rceil = 2$, and $k \ge 3$ or $n \ge 4k$, 
then by a theorem of Fujita et al.~\cite{FMTY2006}, 
$G$ contains $k$ disjoint cycles. 
Thus, 
we may assume that 
$k = 2$, $\lceil m / k \rceil = 2$ and $n = 4k-1 = 7$. 
Then, 
$\delta(G) \ge m \ge k + 1 = 3 = 2k - 1$ 
and 
$\sigma_{3}(G) > \frac{3}{2}(n-1) = 6k-3 = 9$. 
Since $n = 7$ and $\sigma_{3}(G) > 9$, 
it follows that $\alpha (G) \le 3 = n - 2k$ 
and 
$G$ is not a wheel. 
Hence, 
by Theorem~\ref{Thm:KKY2017}, 
$G$ contains two disjoint cycles. 
Thus, the lemma follows. 
\qed

Recall that $\sigma_{t}^{m}(G) \ge \frac{t}{m} \cdot \sigma_{m}(G)$ for $m \ge t \ge 2$, 
and hence 
Theorem~\ref{thm:Yamashita-type for 2-factor with k cycles (Partition)} and Lemma~\ref{lem:cycle packing II} 
lead to the following.

\begin{cor}
\label{cor:Bondy-type for 2-factor with k cycles}
Let $k$ be a positive integer, 
and let $G$ be an $m$-connected graph of order $n \ge 4k-1$. 
If $\sigma_{\lceil m/k \rceil + 1}(G) > \frac{1}{2}(\lceil m/k \rceil + 1) (n-1)$, 
then $G$ has a 2-factor with exactly $k$ cycles. 
\end{cor}




\begin{thebibliography}{99}

\bibitem{Ainouche1992}
A.~Ainouche,
\textit{An improvement of Fraisse's sufficient condition for hamiltonian graphs},
J.~Graph Theory \textbf{16} (1992) 529--543.


\bibitem{Bondy1978} 
J.A.~Bondy,
\textit{A remark on two sufficient conditions for hamilton cycles},
Discrete Math. \textbf{22} (1978) 191--193.


\bibitem{Bondy1980}
J.A.~Bondy,
\textit{Longest paths and cycles in graphs with high degree},
Research Report CORR 80-16, Department of Combinatorics and Optimization, 
University of Waterloo, Waterloo, Ontario, Canada (1980).


\bibitem{BMBook2008}
J.A.~Bondy, U.S.R.~Murty, 
``GraphTheory,'' Springer-Verlag, London, 2008.


\bibitem{BCFGL1997}
S.~Brandt, G.~Chen, R.~Faudree, R.J.~Gould, L.~Lesniak, 
\textit{Degree conditions for 2-factors}, 
J.~Graph Theory \textbf{24} (1997) 165--173.


\bibitem {CGKOSS2007}
G.~Chen, R.J.~Gould, K.~Kawarabayashi, K.~Ota, A.~Saito, I.~Schiermeyer, 
\textit{The Chv\'{a}tal-Erd\H{o}s condition and 2-factors with a specified number of components},
Discuss. Math. Graph Theory \textbf{27} (2007) 401--407. 


\bibitem {CSS2013}
G.~Chen, A.~Saito, S.~Shan, 
\textit{The existence of a 2-factor in a graph satisfying the local Chv\'{a}tal-Erd\H{o}s condition},
SIAM J. Discrete Math. \textbf{27} (2013) 1788--1799. 


\bibitem {CE1972}
V.~Chv\'{a}tal, P.~Erd\H{o}s, 
\textit{A note on hamiltonian circuits},
Discrete Math. \textbf{2} (1972) 111--113.




\bibitem{Enomoto1998}
H.~Enomoto, 
\textit{On the existence of disjoint cycles in a graph},
Combinatorica \textbf{18} (1998) 487--492.


\bibitem{FMTY2006}
S.~Fujita, H.~Matsumura, M.~Tsugaki, T.~Yamashita, 
\textit{Degree sum conditions and vertex-disjoint cycles in a graph},
Australas. J. Combin \textbf{35} (2006) 237--251.


\bibitem{GSurvey2003}
R.J.~Gould, 
\textit{Advances on the hamiltonian problem -- a survey}, 
Graphs and Combin. \textbf{19} (2003) 7--52. 


\bibitem{KY2003} 
A.~Kaneko, K.~Yoshimoto, 
\textit{A 2-factor with two components of a graph satisfying the Chv\'{a}tal-Erd\H{o}s condition}, 
J.~Graph Theory \textbf{43} (2003) 269--279.


\bibitem{KKY2017}
H.A.~Kierstead, A.V.~Kostochka, E.C.~Yeager, 
\textit{On the Corr\'{a}di-Hajnal theorem and a question of Dirac}, 
J. Combin. Theory Ser. B \textbf{122} (2017) 121--148. 


\bibitem{LSurvey2013}
H.~Li, 
\textit{Generalizations of Dirac's theorem in hamiltonian graph theory -- a survey}, 
Discrete Math.~\textbf{313} (2013) 2034--2053.


\bibitem{Ore1960}
O.~Ore,
\textit{Note on hamilton circuits},
Amer. Math. Monthly \textbf{67} (1960) 55.


\bibitem{Wang1999}
H.~Wang,
\textit{On the maximum number of independent cycles in a graph},
Discrete Math.~\textbf{205} (1999) 183--190.


\bibitem{Yamashita2008}
T.~Yamashita,
\textit{On degree sum conditions for long cycles and cycles through specified vertices},
Discrete Math.~\textbf{308} (2008) 6584--6587.

\end{thebibliography}
\end{document}